# ARCHIPELAGO GROUPS

GREGORY R. CONNER, WOLFRAM HOJKA, AND MARK MEILSTRUP

ABSTRACT. The classical *archipelago* is a non-contractible subset of $\mathbb{R}^3$ which is homeomorphic to a disk except at one non-manifold point. Its fundamental group, $\mathscr{A}$, is the quotient of the *topologist's product* of $\mathbb{Z}$, the fundamental group of the shrinking wedge of countably many copies of the circle (the Hawaiian earring), modulo the corresponding free product. We show $\mathscr{A}$ is locally free, not indicable, and has the rationals both as a subgroup and a quotient group. Replacing $\mathbb{Z}$ with arbitrary groups yields the notion of *archipelago groups*.

Surprisingly, every archipelago of countable groups is isomorphic to either $\mathscr{A}(\mathbb{Z})$ or $\mathscr{A}(\mathbb{Z}_2)$, the cases where the archipelago is built from circles or projective planes respectively. We conjecture that these two groups are isomorphic and prove that for large enough cardinalities of $G_i$, $\mathscr{A}(G_i)$ is not isomorphic to either.

## 1. INTRODUCTION

In algebraic topology one of the prototypical two-dimensional spaces eliciting unusual properties in its fundamental group is the harmonic archipelago: commonly seen as a shrinking bouquet of circles, with discs glued in between consecutive circles that bulge out to a constant height. Perhaps it should come as no surprise that this space has interesting features, as it is homeomorphic to the *reduced metric suspension* of the graph of the topologist's sine curve $y = \sin(1/x)$, a troublemaker well-known from analysis and topology (Proposition 12).

For our discussions it will be most useful to switch to a homotopy equivalent construction, namely to the mapping cone of the natural continuous map from a regular wedge of spaces to the same wedge but with a strong topology. We will call this mapping cone $A$ the archipelago space for the sequence of spaces (Definition 4) and its fundamental group then turns out to be the penultimate term in the exact sequence

$$1 \to \langle\!\langle \pi_1(\bigvee_n X_n)\rangle\!\rangle \to \pi_1\big(\bigotimes_n X_n\big) \to \pi_1(A) \to 1,$$

where "$\langle\!\langle\ \rangle\!\rangle$" denotes taking the normal closure in the next term of the sequence.

If the spaces $X_n$ are all taken to be circles, the resulting mapping cone is homotopy equivalent to the harmonic archipelago. As another example, consider the archipelago space of a sequence of projective planes. Our main theorem will show that, amazingly, the fundamental group of *any* well-behaved archipelago space, e.g. where each $X_n$ is a locally finite CW-complex, is either trivial or isomorphic to one of these two examples.

If the spaces are locally connected and first-countable at the basepoint (or alternatively, if the wedges are understood as *homotopy* colimits instead of regular ones), the first nontrivial group in the above sequence is isomorphic to the normal closure of the ordinary free product $\ast_n \pi_1 X_n$; the second one is the *topologist's product* $\otimes_n \pi_1 X_n$, first studied

2010 *Mathematics Subject Classification.* Primary 55Q20, 20E06; Secondary 57M30, 57M05, 20F05.

This work was supported by the Simons Foundation Grant 246221 and by the Austrian Science Foundation FWF project S9612.

by Higman [Hig52, Section 6] and Griffiths [Gri56, Section 4] in the 50's (and is also known as $\sigma$-product or free complete product). This is an analogue of the finite concatenations in a free product to infinite concatenations. It can be intuitively understood by using the concept of infinite words, see [Eda92, Definition 1.1]. Also, the free product is contained in it as a subgroup, just as the direct sum embeds in the direct product for abelian groups.

It is possible to define an *archipelago group* $\mathscr{A}(G_n)$ of a sequence of groups in purely algebraic terms as the quotient $\bigoplus_n G_n / \langle\!\langle *_n G_n \rangle\!\rangle$ (as in Definition 6), and in the special case of all groups being $\mathbb{Z}$, this is the fundamental group of the harmonic archipelago.

These groups allow the following remarkable classification.

**Theorem A.** *Let $(G_n)_{n \in \mathbb{N}}$ be a collection of nontrivial countable (possibly finite) groups. If only finitely many of the groups $G_n$ have elements of order 2, then*

$$\mathscr{A}(G_n) \simeq \mathscr{A}(\mathbb{Z}).$$

*If infinitely many of the groups $G_n$ have elements of order 2, then*

$$\mathscr{A}(G_n) \simeq \mathscr{A}(\mathbb{Z}_2).$$

That is, any such group is isomorphic to either of two prototypes! In particular, the groups $\mathscr{A}(\mathbb{Z})$, $\mathscr{A}(\mathbb{Q})$, and $\mathscr{A}(\mathbb{Z}_3)$ are all isomorphic. The next theorem collects various properties about this standard case.

**Theorem B.** *The fundamental group $\mathscr{A}(\mathbb{Z})$ of the harmonic archipelago has the following properties.*

(1) *The group $\mathbb{Q}$ of rational numbers embeds as a subgroup in $\mathscr{A}(\mathbb{Z})$;*
(2) *$\mathscr{A}(\mathbb{Z})$ does not have an infinite cyclic quotient (i.e. is not indicable); and*
(3) *$\mathscr{A}(\mathbb{Z})$ is locally free.*

The proof of Theorem A is astonishing in that it only employs arbitrary set theoretic bijections between groups to derive an *isomorphism* between the archipelago groups. This "trick" fails for involution elements, leading to the two possible cases above. To make matters even more interesting, we will see later that the group $\mathscr{A}(\mathbb{Z}_2)$ is also torsion-free, that $\mathscr{A}(\mathbb{Z})$ embeds in $\mathscr{A}(\mathbb{Z}_2)$ so consequently $\mathscr{A}(\mathbb{Z}_2)$ contains a copy of the rationals, and that it is also not indicable. Furthermore in [HH13] it is shown that $\mathscr{A}(\mathbb{Z}_2)$ is locally free, naturally leading one to ask if the two groups are the same:

**Question 1.** *Is $\mathscr{A}(\mathbb{Z}) \simeq \mathscr{A}(\mathbb{Z}_2)$?*

Indeed, more generally, any archipelago group is completely encoded by a sequence of cardinals (Theorem 10); this bears a certain superficial resemblance to the classical situation for abelian groups, as [Hul62], together with standard results in [Kap54], shows the quotient $\prod_n G_n / \bigoplus_n G_n$ to also be depending only on a sequence of cardinal numbers, albeit a quite different one (see the extended discussion after Theorem 10 in section 3).

However, distinct cardinalities do not necessarily lead to different archipelago groups, and it is also conceivable that the following two groups are isomorphic.

**Question 2.** *Is $\mathscr{A}(\mathbb{Z}) \simeq \mathscr{A}(\mathbb{R})$?*

## 2. THE TOPOLOGICAL VIEWPOINT

For a collection of spaces $X_n$ with basepoints $p_n$ the wedge $\bigvee_n X_n$ is their coproduct in the category of pointed spaces. As a set, it can be naturally considered as the subset of the product $\prod_n X_n$ consisting of all points $(a_n)_n$ such that $a_n = p_n$ for all but at most one $n$, the "main axes". If the product is given the box topology this yields the usual weak wedge topology. If however the product is given the standard (Tychonoff) topology, that same subset will be called the *shrinking wedge* $\varoslash_n X_n$.

This information can be encoded in the left half of following diagram:

$$\begin{array}{ccc} \bigvee_n X_n \hookrightarrow \prod_n^{\text{BOX}} X_n & & *_n G_n \twoheadrightarrow \sum_n G_n \\ f \downarrow \quad \downarrow & \xrightarrow{\pi_1} & \uparrow \quad \uparrow \\ \varoslash_n X_n \hookrightarrow \prod_n X_n & & \circledast_n G_n \twoheadrightarrow \prod_n G_n \end{array} \quad (\#)$$

The horizontal maps are the embeddings, the vertical maps are continuous bijections. Assuming the spaces are *nice* at the basepoint (see below), setting $G_n := \pi_1(X_n)$ and applying the functor $\pi_1$ yields the right part of the diagram together with the induced maps between the fundamental groups. On this side, the vertical maps are embeddings, the horizontal ones are onto. The symbol $\sum$ retrieves the set of elements in the product with finite support. For abelian groups we will use the more common $\bigoplus$, denoting the same object.

As an example, recall that the Hawaiian earring is the planar set consisting of circles $c_n$ of radius $1/n$ centered at $(0, 1/n)$. Note that each circle contains the origin, and equivalently this space is the one-point compactification of a sequence of open arcs. It is also a shrinking wedge $\varoslash_{n \in \mathbb{N}} S^1$ of infinitely many circles and its fundamental group turns out to be $\circledast_{n \in \mathbb{N}} \mathbb{Z}$, each copy of $\mathbb{Z}$ corresponding to one circle of the earring.

Algebraically, the most straight-forward method of defining the topologist's product for countable families is simply the following:

$$\circledast_i G_i := \bigcap_i \left( G_i * \varprojlim_n *_{1 \leq j \leq n, \, j \neq i} G_j \right).$$

This group was first implicitly conceived by Higman in [Hig52, Section 6] while studying subgroups of $\varprojlim_n *_{1 \leq j \leq n} G_j$, a group he called the *unrestricted free product* of the $G_i$. Griffiths then showed in [Gri56, Theorem 6.3] (with a corrected proof due to Morgan and Morrison in [MM86, p. 563, Infinite van Kampen Theorem]) the relation

$$\pi_1(\varoslash_i X_i) = \circledast_i \pi_1(X_i)$$

between the topologist's product and the shrinking wedge of spaces that satisfy some local properties: being locally simply connected and first countable at each basepoint.

Contrast this with the standard fact that

$$\pi_1(\bigvee_i X_i) = *_i \pi_1(X_i)$$

holds under the same assumptions on the $X_i$.

These local requirements can be avoided by a simple procedure. To each space $(X_i, p_i)$ attach an arc to $p_i$ and shift the basepoint to the arc's other end. Let $(\tilde{X}_i, \tilde{p}_i)$ denote the thus modified space, then we can define a *homotopy shrinking wedge* $\varobar^{\mathrm{H}}_i(X_i, p_i) := \varobar(\tilde{X}_i, \tilde{p}_i)$. Now, with this notation,

$$\pi_1(\varobar^{\mathrm{H}}_i X_i) = \circledast_i \pi_1(X_i)$$

holds for arbitrary spaces $X_i$.

Similarly, one can define a *homotopy wedge* by the same method. Note, that the difference between a regular wedge and a homotopy wedge is reflected in the distinction between taking the colimit or the *homotopy* colimit of the diagram of a wedge, i.e. the collection of maps $p_i : \{x\} \to X_i$ that select the basepoint in each space. Namely, following e.g. the construction in [Vog73, Definition 1.1], the homotopy colimit of this simple diagram will have a copy of the unit interval $I$ for each map $p_i$, with one endpoint identified with $x$, the other with the point $p_i(x)$. The abstract theory of homotopy colimits would then also yield the above isomorphism $\pi_1(\bigvee_i X_i) = *_i \pi_1(X_i)$, as the free product is the (homotopy) colimit of the induced diagram of groups. This is a special case of Theorem 1.1 in [Far04], where these concepts are discussed in more generality.

Another interesting property of the topologist's product is that it can be interpreted as an infinite word structure.

**Definition 3.** For a family $(G_n)_{n \in \mathbb{N}}$ of groups, an *infinite word* is a map $w : L \to \bigcup G_n \smallsetminus \{1\}$ from a countable linearly ordered set $L$ to the disjoint union of the non-identity elements of the $G_n$ where the preimage of every $G_n$ is a finite set. Multiplication is simply concatenation, and inverses are given by inverting the order of the word and replacing each element by its inverse in $G_n$.

Corresponding to the case of finite words, there is a natural notion of cancellation within an infinite word (see [CC00, Definition 3.4] and [Eda92] for more background on this concept). Cancellation induces an equivalence relation on the set of infinite words, and the classes together with the above operations form the *topologist's product* $\circledast_n G_n$.

We return our attention to the harmonic archipelago, as defined by Bogley and Sieradski [BS98]. This space can be depicted as an earring with discs $D_i$ attached along each boundary $c_i c_{i+1}^{-1}$ where an interior point of $D_i$ is raised to height $z = 1$ in $\mathbb{R}^3$, see Figure 1.

Recall that for a map $f : X \to Y$, the mapping cone $C_f$ is defined as $X \times I \cup Y$ with every $(x, 0)$ identified with $f(x)$ and $X \times \{1\}$ collapsed to a point. Then the archipelago is homotopy equivalent to the mapping cone $C_f$ of the continuous bijective map $f : \bigvee_n S^1 \to \varobar_n S^1$ from the wedge to the shrinking wedge of circles, corresponding to the leftmost map in (#), see Proposition 13 for a detailed argument.

This allows us to generalize the notion to arbitrary spaces.

**Definition 4.** The *archipelago space* $A$ of a family of spaces $X_n$ is the mapping cone $C_f$ for $f : \bigvee_n X_n \to \varobar_n X_n$.

**Theorem 5.** *The fundamental group of the archipelago space* $A = C_f$ *is determined by the fundamental groups* $G_n := \pi_1(X_n)$ *as the cokernel of* $f_*$ *in the sequence*

$$\pi_1(\bigvee_n X_n) = *_n G_n \xhookrightarrow{f_*} \pi_1(\varobar_n X_n) = \circledast_n G_n \twoheadrightarrow \pi_1(A).$$

*It can directly be expressed as the quotient*

$$\pi_1(A) = \circledast_i G_i / \langle\!\langle *_i G_i \rangle\!\rangle.$$

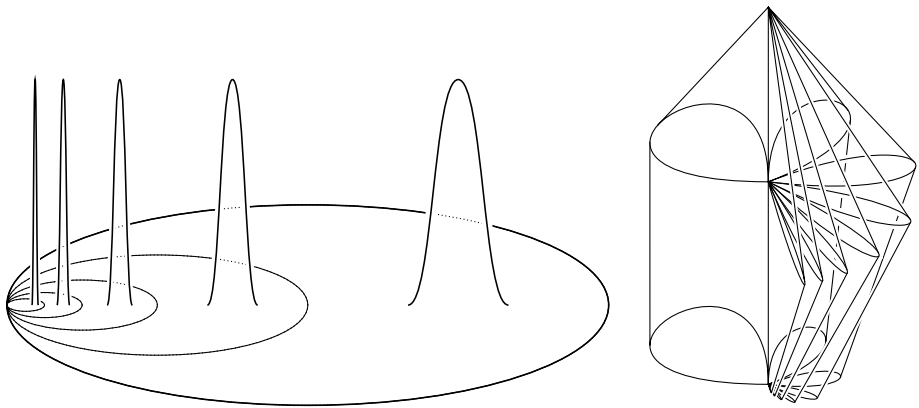

Figure 1: (left) The harmonic archipelago. (right) A homotopy equivalent realization as the mapping cone of $f : \bigvee_n S^1 \to \obigtimes_n S^1$.

The proof of the theorem is a simple application of van-Kampen's theorem. The homomorphism $f_*$ induced by $f$ coincides precisely with the canonical inclusion map of the free product in the topologist's product, from which the alternative quotient representation of the fundamental group above follows.

In the introduction we have mentioned how this quotient is interesting in purely algebraic terms. The preceding theorem motivates its name:

**Definition 6.** Given a countable collection of groups $G_n$, the *archipelago group* $\mathscr{A}(G_n)$ is given by the topologist's product of the $G_n$'s, modulo their free product:
$$\mathscr{A}(G_n) := \obigtimes_n G_n / \langle\!\langle *_n G_n \rangle\!\rangle.$$
If each $G_n$ is the same group $G$, we write $\mathscr{A}(G)$.

## 3. Further Results

Theorem B can be extended in various ways. Finding a subgroup isomorphic to $\mathbb{Q}$ is similar in spirit to [BZ12] where the same is shown for Griffiths' double cone space. It is also a corollary to the theorem in [Hoj13] that *every* countable locally free group embeds in the archipelago.

The second property of Theorem B is not unique to $\mathscr{A}(\mathbb{Z})$, but holds for any archipelago group.

**Proposition 7.** *For any countable collection $(G_n)_{n \in \mathbb{N}}$ of groups, $\mathscr{A}(G_n)$ does not have an infinite cyclic quotient. In other words,* $\mathrm{Hom}(\mathscr{A}(G_n), \mathbb{Z})$ *is trivial.*

Similarly, neither is the third property unique, every archipelago group is locally free according to [HH13]. As an intermediate step, here we will be satisfied with the following weaker fact.

**Proposition 8.** *The group $\mathscr{A}(G_n)$ is torsion-free.*

Question 1 left as a possibility that every archipelago group composed of countable groups $G_n$ is isomorphic to that of the standard harmonic archipelago. But certainly, not *all* archipelago groups are equal; for reasons of cardinality there is an infinite class of them, all distinct:

**Theorem 9.** *Suppose all groups $G_n$ are of the same cardinality $\kappa$, then the archipelago group $\mathscr{A}(G_n)$ has cardinality $\kappa^{\aleph_0}$. Hence there is more than a set's worth of distinct archipelago groups. In particular, the archipelago groups obtained from free groups of $\beth_{\alpha+1}$ many generators are all distinct for $\alpha \geq 0$.*

This however still allows for the two groups from Question 2, $\mathscr{A}(\mathbb{Z})$ and $\mathscr{A}(\mathbb{R})$, of cardinality $\leq \mathfrak{c} = \beth_1$ to be isomorphic to each other.

On the other hand there is the following positive result on representing archipelago groups by cardinal numbers.

**Theorem 10.** *Let $\kappa_n$ denote the cardinality of $G_n \smallsetminus \{1\}$, and $\lambda$ that of the set of indices $i$ such that $G_i$ contains an involution. Then $\mathscr{A}(G_n)$ is determined by these countably many cardinal numbers in the following way: if $\lambda < \aleph_0$,*

$$\mathscr{A}(G_n) \simeq \mathscr{A}(*_{\kappa_n} \mathbb{Z})$$

*and if $\lambda = \aleph_0$,*

$$\mathscr{A}(G_n) \simeq \mathscr{A}(*_{\kappa_n} \mathbb{Z}_2).$$

Theorem 10 compares to classical results for abelian groups: In [Hul62] Hulanicki shows that the quotient $\prod_i G_i / \bigoplus_i G_i$ is algebraically compact, a slightly generalized statement is in [Fuc63]. Using standard techniques due to Kaplansky [Kap54], the quotient then has a representation as

$$\prod_i G_i / \bigoplus_i G_i \simeq \bigoplus_{\mathfrak{k}} \mathbb{Q} \oplus \bigoplus_{p \text{ prime}} \bigoplus_{\mathfrak{l}_p} \mathbb{Z}(p^\infty) \oplus \prod_{p \text{ prime}} A_p$$

where

$$A_p \simeq p\text{-adic completion of } \left( \bigoplus_{\mathfrak{m}_{p,0}} \mathbb{J}_p \oplus \bigoplus_{k=1}^{\infty} \bigoplus_{\mathfrak{m}_{p,k}} \mathbb{Z}/p^k\mathbb{Z} \right)$$

that only depends on the countably many cardinals $\mathfrak{k}$, $\mathfrak{l}_p$, $\mathfrak{m}_{p,k}$ (for $p$ prime, $k \in \mathbb{N}$), with $\mathbb{Z}(p^\infty)$ the $p$-quasicyclic group and $\mathbb{J}_p$ the $p$-adic completion of $\mathbb{Z}$. The relevant arguments are the statements 21.3, 23.1, 40.1, and the remark after 40.2 in Fuchs' book [Fuc70] on that subject.

Similarly, the nonabelian quotient $\mathscr{A}(G_n) = \circledast_n G_n / \langle\!\langle *_n G_n \rangle\!\rangle$ depends only on the cardinalities of the groups and the presence of 2-torsion elements. However there is no direct correspondence between the two sets of cardinal numbers.

As $\circledast_n G_n$ maps onto $\prod_n G_n$, one might wonder if for abelian groups the canonical map $\circledast_n G_n / \langle\!\langle *_n G_n \rangle\!\rangle \to \prod_n G_n / \bigoplus_n G_n$ between the quotients is simply induced by abelianization of the group – something that is suggested by the case of $\mathscr{A}(\mathbb{Z})$: here the abelianization coincides with the first singular homology group of the harmonic archipelago. This homology group was e.g. investigated in [KR12] where its torsion-freeness is used. They reference a personal note by Eda, but that property also follows from the stronger result of our Theorem B, part (3). This group is isomorphic to $\prod_n \mathbb{Z} / \bigoplus_n \mathbb{Z}$. Yet the abelianization map turns out to be more complicated than the canonical map above.

The issue is more readily seen for e.g. finite cyclic groups, as these give rise to torsion in the abelian quotient. To wit, let $\mathbb{Z}_k$ denote the cyclic group of order $k$, then $(1,1,1,\ldots) \in$ $\prod_n \mathbb{Z}_k$ gives an element of order $k$ in $\prod_n \mathbb{Z}_k / \bigoplus_n \mathbb{Z}_k$. But if the order $k$ is coprime to 2, the nonabelian quotient $\circledast_n \mathbb{Z}_k / \langle\!\langle *_n \mathbb{Z}_k \rangle\!\rangle$ is isomorphic to $\mathscr{A}(\mathbb{Z})$ by Theorem A and thus locally free by Theorem B. Hence its abelianization is torsion-free.

In [HH13] it is shown that the abelianizations of $\mathscr{A}(\mathbb{Z})$, $\mathscr{A}(\mathbb{Z}_2)$, and $\mathscr{A}(\mathbb{R})$ are in fact all equal to each other; perhaps offering some support to Questions 1 and 2.

## 4. Proofs

We begin by elaborating on the two topological characterizations of the archipelago, first as a reduced suspension, then as a mapping cone.

**Definition 11.** The *topologist's sine curve* is the subset of the plane $\Gamma := \{(x,y) : x \in (0, 1/\pi], y = \sin(1/x)\} \cup \{(0,0)\}$.

**Proposition 12.** *The harmonic archipelago is homeomorphic to the reduced metric suspension of the topologist's sine curve.*

*Proof.* Let $A$ be a closed subset of the metric space $X$, then collapsing $A$ yields a natural metric on $X/A$ (not necessarily compatible with the quotient topology) by $d([x],[y]) := \min\{d(x,y), d(x,A) + d(y,A)\}$.

The reduced metric suspension $\Sigma_M \Gamma$ at the basepoint $(0,0)$ is the bijective metric image of the normal reduced suspension $\Sigma\Gamma$ and can be realized by the set $\Gamma \times [-\pi, \pi] \subseteq \mathbb{R}^3$ with the subsets $\Gamma \times \{-\pi, \pi\}$ and $\{(0,0)\} \times [-\pi, \pi]$ collapsed to a point. The *metric quotient* map $q : \Gamma \times [-\pi, \pi] \to \Sigma_M \Gamma$ can be realized as the continuous function taking $(0, 0, r)$ to the origin and for $x \neq 0$

$$q : (x, y, r) \mapsto \left(x \cdot (1 + \cos r), x \cdot \sin r, \cos(r/2) \cdot \underbrace{|\sin(1/x)|}_{y}\right),$$

whose image in $\mathbb{R}^3$ is a harmonic archipelago. $\square$

*Remark.* The reduced metric suspension of a metric space is homotopy equivalent to its reduced suspension with the quotient topology. Hence, the harmonic archipelago is homotopy equivalent to the reduced suspension of the topologist's sine curve.

**Proposition 13.** *The harmonic archipelago is homotopy equivalent to the mapping cone of the canonical map from the wedge of circles to their shrinking wedge.*

*Proof.* We will approach this by interpreting both spaces as adjunctions of 2-cells to an Hawaiian earring, i.e. as *relative* CW-complexes.

First notice, that there are two distinct ways to define a topology on the harmonic archipelago: the humps between the circles in the base Hawaiian earring can be either given the topology from the metric embedding into $\mathbb{R}^3$ or a weak topology as an adjunction space with attached discs. The base earring in the metric archipelago allows a neighbourhood deformation retract – by contracting an annulus in each disc toward its bounding circle in the base – inducing a homotopy equivalence between the two variants.

Adopting a more general viewpoint, this weak archipelago is the pushout of a union of discs $\bigsqcup D^2$ and an earring $\bigodot S^1$ along the adjunction map $k : \bigsqcup S^1 \to \bigodot S^1$ that takes the $i$-th circle $c_i$ to $b_i \cdot b_{i+1}^{-1} \cdot 1 \cdot 1$ (cf. the front square of the diagram in Figure 2).

On the other hand, the mapping cone may be thought of as a pushout as well, one of a union of discs and an earring *with an arc attached*, $I \vee \bigodot S^1$; this time with $j : \bigsqcup S^1 \to$ $I \vee \bigvee S^1$ mapping the $i$-th circle $c_i$ to the path $a_i \cdot h \cdot h^{-1}$ (with $h$ the path running up the line segment $I$).

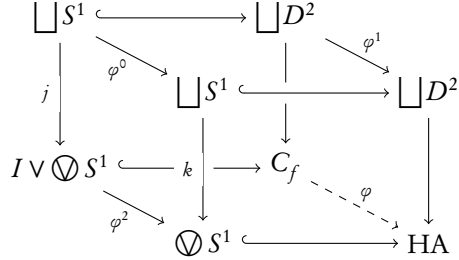
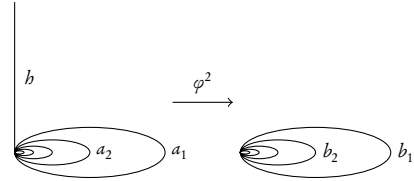

Figure 2: The pushout diagrams for the gluing theorem.

Figure 3: Domain and codomain of the homotopy equivalence $\varphi^2$.

This second pushout corresponds to the back square in the diagram. Now we can add homotopy equivalences $\varphi^0, \varphi^1, \varphi^2$, the first two simply being identity maps, the last one, $\varphi^2$, being defined as collapsing the line segment of $h$ to the basepoint and mapping the circle $a_i$ to the path $b_i \cdot b_{i+1}^{-1}$ (see Figure 3). As $k \circ \varphi^0(c_i) = b_i \cdot b_{i+1}^{-1} \cdot 1 \cdot 1 = \varphi^2(a_i \cdot h \cdot h^{-1}) = \varphi^2 \circ j(c_i)$, the diagram is commutative. Thus the gluing theorem [Bro06, 7.5.7] shows that the canonical map $\varphi$ is a homotopy equivalence between the mapping cone and the weak archipelago. Combining this map with the equivalence above to the (metric) harmonic archipelago yields the desired result. □

A helpful tool in deriving algebraic properties for an archipelago is representing the group as a direct limit.

**Proposition 14.** *The archipelago group $\mathscr{A}(G_n)$ is a direct limit,*

$$\mathscr{A}(G_n) = \varinjlim_k \circledast_{i \geq k} G_i,$$

*where the bonding maps are the quotient maps*

$$\circledast_{i \geq k} G_i \to \left( \circledast_{i \geq k} G_i \right) / \langle\!\langle G_k \rangle\!\rangle \simeq \circledast_{i \geq k+1} G_i.$$

*Proof.* Note that any element of the free product $*_{i \geq 1} G_i$ becomes trivial after the application of finitely many bonding maps. Thus any group mapped to from this directed system of groups must be a factor of $\circledast_n G_n / \langle\!\langle *_n G_n \rangle\!\rangle = \mathscr{A}(G_n)$, which then satisfies the universal property of the direct limit. □

In the following lemma, we will use functions between groups that have some of the properties of homomorphisms, but are somewhat less restrictive. We call a function $\varphi : G \to H$ *identity-preserving* if $\varphi(e_G) = e_H$. Similarly, $\varphi$ is *inverse-preserving* if for all $g \in G$, $\varphi(g^{-1}) = (\varphi(g))^{-1}$.

**Lemma 15.** *If for every $i$ there exists an identity-preserving and inverse-preserving function $\varphi_i : G_i \to H_i$ (not necessarily a homomorphism), then the functions $\varphi_i$ induce a homomorphism $\Phi : \mathscr{A}(G_i) \to \mathscr{A}(H_i)$.*

*Remark.* Note that we do not require $\varphi_i$ to preserve multiplication.

*Remark.* For the technical details of the proof it is necessary to recapitulate certain properties of the word calculus used to describe the topologist's product. Consider a word $w : L \to \bigcup G_i \smallsetminus \{1\}$ as in Definition 3 and write the concatenation of words $u$ and $v$ as $u \cdot v$. As the preimage of every group $G_i$ is finite, there is a projection $p_n(w)$ of $w$ to a finite word in $*_{i=1}^n G_i$ for each $n$. A word $w$ is called *reduced* if

(i) for each nontrivial subword $w'$ (i.e. $w = u \cdot w' \cdot v$) there exists some projection $p_n(w')$ that is nontrivial; and
(ii) no two consecutive points in the order $L$ are mapped to $g_1, g_2$ in the same group $G_i$.

In fact, every element of the topologist's product (seen as a class of words) can be represented by a reduced word, unique up to order isomorphism, and one can naturally identify the elements of the group with the reduced words, cf. Definition 3 and [Eda92, Theorem 1.4]. Let $R(w)$ denote this reduced representative of a word $w$.

Of particular interest will be finding $R(u \cdot v)$ for already reduced $u$ and $v$, which can be accomplished by two simple reduction steps. It is possible to write the two reduced words as concatenations $u = a \cdot g_1 \cdot x$ and $v = y \cdot g_2 \cdot b$ satisfying: firstly, $x$ is the maximal (possibly trivial) terminal subword of $u$ with its inverse word $y = x^{-1}$ as an initial subword of $v$; secondly the letters $g_1, g_2$ are either both trivial, or they are elements of the same factor $G_i$ whose product $g_1 g_2$ is nontrivial. Then $R(u \cdot v) = a \cdot (g_1 g_2) \cdot b$ is reduced as a concatenation.

*Proof of the Lemma.* The union of the functions $\varphi_i$ gives a set function $\varphi_0$ between the disjoint unions $\bigcup_i G_i$ and $\bigcup_i H_i$. This function $\varphi_0$ will induce a map (not necessarily a homomorphism) $\varphi$ between the topologist's products, which will in turn induce our desired homomorphism $\Phi : \mathscr{A}(G_i) \to \mathscr{A}(H_i)$, as we describe below.

The set function $\varphi_W$ is a letter replacement defined by $\varphi_0$, in the following manner: an infinite word in the $G_n$ is mapped to an infinite word in the $H_n$ of the same order type, by taking each letter $a \in G_j$ to $\varphi_0(a) = \varphi_j(a)$. For example, $a_1 a_2 a_3 \ldots \mapsto \varphi_0(a_1)\varphi_0(a_2)\varphi_0(a_3)\ldots$. Generally, an infinite word $w : L \to \bigcup G_i \smallsetminus \{1\}$ is mapped to

$$\varphi_W(w) := \varphi_0 \circ w \;:\; L \to \bigcup H_n \smallsetminus \{1\},$$

an infinite word in the $H_n$.

By considering group elements as reduced words, $\varphi_W$ also induces a map on the topologist's product. Define $\varphi : \circledast_n G_n \to \circledast_n H_n$, $\varphi(u) := R \circ \varphi_W(u)$.

This $\varphi$ will not be a homomorphism; it does however resemble one, at least for reduced concatenations. If a product $uv$ of group elements is equal to its concatenation $u \cdot v$ (or equivalently, $R(u \cdot v) = u \cdot v$), then

$$\varphi(uv) = \varphi(u)\varphi(v). \tag{\dag}$$

This follows from the chain of equalities $\varphi(uv) = R \circ \varphi_W(u \cdot v) = R(\varphi_W(u) \cdot \varphi_W(v)) = R(\varphi(u) \cdot \varphi(v)) = \varphi(u)\varphi(v)$.

We claim that $\varphi$ induces a homomorphism $\Phi$ on the induced archipelago groups:

$$\begin{array}{ccc} \circledast_n G_n & \xrightarrow{\varphi} & \circledast_n H_n \\ \downarrow & & \downarrow \\ \mathscr{A}(G_n) & \xrightarrow{\Phi} & \mathscr{A}(H_n) \end{array} \tag{\ddag}$$

We will write $w \sim v$ to denote the equivalence between elements in the topologist's product projecting down to the same element in the archipelago group (i.e. in the quotient by the free product). It needs to be shown that $\Phi$ is well-defined, in other words, that $w \sim v$ implies $\varphi(w) \sim \varphi(v)$; and that it is actually a homomorphism.

**Claim.** *For reduced words $u$ and $v$, $\varphi(R(u \cdot v)) \sim \varphi(u)\varphi(v)$.*

By the remark it is possible to write the reduced words as concatenations $u = a \cdot g_1 \cdot x$ and $v = y \cdot g_2 \cdot b$, with $y = x^{-1}$ the inverse word of $x$ and $g_1, g_2$ in some $G_i$. Then $R(u \cdot v) = a \cdot (g_1 g_2) \cdot b$ is reduced as a concatenation of reduced words, and thus from (†),

$$\varphi(R(u \cdot v)) = \varphi(a(g_1 g_2)b) = \varphi(a)\varphi(g_1 g_2)\varphi(b).$$

Further note, $h \sim 1$ for $h \in H_i$ and $\varphi$ preserves inverses, thus

$$\begin{aligned}
\varphi(R(u \cdot v)) &= \varphi(a)\varphi(g_1 g_2)\varphi(b) \\
&\sim \varphi(a)\varphi(b) \\
&\sim \varphi(a)\varphi(g_1)\varphi(g_2)\varphi(b) \\
&= \varphi(a)\varphi(g_1)\varphi(x)\varphi(x)^{-1}\varphi(g_2)\varphi(b) \\
&= \varphi(a)\varphi(g_1)\varphi(x)\varphi(x^{-1})\varphi(g_2)\varphi(b) \\
&= \varphi(a g_1 x)\varphi(y g_2 b) \\
&= \varphi(u)\varphi(v),
\end{aligned}$$

as claimed.

Now suppose $u \sim v$ for two elements $u, v \in \circledast_n G_n$. Then using the direct limit representation of the archipelago there exists an index $j$ so that $\tau_j(u) = \tau_j(v)$ for the canonical projection $\tau_j : \circledast_n G_n \to \circledast_{n>j} G_n$. By virtue of the free decomposition

$$\circledast_n G_n = (G_1 * \ldots * G_j) * \circledast_{n>j} G_n$$

$u$ can be written as a finite product $u_1 \ldots u_r$, with each $u_i$ either in $G_1 * \ldots * G_j$ or $\circledast_{n>j} G_n$. Let $U_1, \ldots, U_k$ denote the subsequence of factors lying in the latter, then their product $U_1 \ldots U_k$ yields precisely $\tau_j(u)$. By recursively applying the claim,

$$\varphi(U_1)\ldots\varphi(U_k) \sim \varphi\bigl(R(U_1 \cdot U_2)\bigr)\varphi(U_3)\ldots\varphi(U_k) \sim \ldots$$
$$\sim \varphi\bigl(R(U_1 \cdot \ldots \cdot U_k)\bigr),$$

and hence

$$\begin{aligned}
\varphi(u) &= \varphi(u_1 \ldots u_r) \\
&= \varphi(u_1)\ldots\varphi(u_r) \\
&\sim \varphi(U_1)\ldots\varphi(U_k) \\
&\sim \varphi\bigl(R(U_1 \cdot \ldots \cdot U_k)\bigr) \\
&= \varphi\bigl(\tau_j(u)\bigr).
\end{aligned}$$

Repeating the same process for $v$ we arrive at $\varphi(u) \sim \varphi\bigl(\tau_j(u)\bigr) = \varphi\bigl(\tau_j(v)\bigr) \sim \varphi(v)$, as desired. So $\Phi$ is well-defined.

It remains to show that $\Phi$ is a homomorphism, i.e.

$$\Phi([u][v]) = \Phi([u])\Phi([v])$$

for $u, v \in \circledast_n G_n$. But for this it suffices that $\varphi(uv) = \varphi(R(u \cdot v)) \sim \varphi(u)\varphi(v)$ holds with respect to their reduced word representations, precisely as stated in the claim. □

**Corollary 16.** *If the maps $\varphi_i$ are bijections, then $\Phi$ is an isomorphism.*

*Proof.* As the same argument now holds mutatis mutandis for the inverses of the original bijections $G_i \to H_i$, the so defined map $\varphi^{-1}$ induces a homomorphism that is clearly inverse to $\Phi$. □

From Corollary 16, we can find many archipelago groups that are isomorphic. The main factors to consider are the cardinalities of the groups involved, and the number of those groups containing 2-torsion. The groups' cardinalities are essential in order to even construct bijections among them. The 2-torsion is the only obstruction to then constructing *inverse-preserving* bijections.

Before turning to the proofs of the theorems themselves, we mention a few basic properties.

**Lemma 17.** *An archipelago group satisfies the following:*

(1) $\mathscr{A}(G_n)$ *is independent of the ordering of the groups, i.e. for any permutation* $f : \mathbb{N} \to \mathbb{N}$, $\mathscr{A}(G_n) \simeq \mathscr{A}(G_{f(n)})$;
(2) $\mathscr{A}(G_n) \simeq \mathscr{A}(G_n)_{n \geq k}$; *and*
(3) $\mathscr{A}(G_n) \simeq \mathscr{A}(G_{2n-1} * G_{2n})$.

*Proof.* First, (1) follows from the fact that both, the topologist's and the free product, share this symmetric property: $\circledast_n G_n \simeq \circledast_n G_{f(n)}$ and $*_n G_n \simeq *_n G_{f(n)}$.

Also, finitely many groups can be split off as a free factor from these products, as $\circledast_n G_n \simeq (*_{n \leq k} G_n) * \circledast_{n > k} G_n$. Thus by cancelling the left factors in

$$\begin{aligned}\mathscr{A}(G_n) = \circledast_n G_n \,/\, \langle\!\langle *_n G_n \rangle\!\rangle &\simeq \\ &\simeq (*_{n<k} G_n * \circledast_{n \geq k} G_n) \,/\, \langle\!\langle *_{n<k} G_n * *_{n \geq k} G_n \rangle\!\rangle \simeq \\ &\simeq \circledast_{n \geq k} G_n \,/\, \langle\!\langle *_{n \geq k} G_n \rangle\!\rangle,\end{aligned}$$

property (2) follows.

Lastly, (3) is inherited from $\circledast_n G_n \simeq \circledast_n (G_{2n-1} * G_{2n})$. □

Indeed, these properties can be combined in the following, slightly more general statements, whose proof is left to the reader.

**Proposition 18.** *Let $A$ be a finite subset of $\mathbb{N}$ and let $(P_n)_{n \in \mathbb{N}}$ be a partition of $\mathbb{N} \setminus A$ into finite sets. Then $\mathscr{A}(G_n) \simeq \mathscr{A}(*_{i \in P_n} G_i)$.*

**Proposition 19.** *Let $(P_n)_{n \in \mathbb{N}}$ be a partition of $\mathbb{N}$ into finite sets. If $\varphi_n : G_n \to *_{i \in P_n} H_i$ is an embedding for all $n \in \mathbb{N}$, then $\mathscr{A}(G_n)$ embeds in $\mathscr{A}(H_n)$.*

The second proposition follows immediately from the first, with the embedding induced as in the diagram (‡). In particular, $\mathscr{A}(\mathbb{Z}_2)$ embeds in $\mathscr{A}(\mathbb{Z})$.

**Lemma 20.** *If $G$ and $H$ are nontrivial, the cardinality of the set of non-involutions of $G * H$ is given by*

$$\left|\{x \in G * H : x^2 \neq 1\}\right| = \max\{\aleph_0, |G|, |H|\}.$$

*If at least one of the two groups contains an involution, then the set of involutions in the product has the same cardinality*

$$\left|\{x \in G * H : x^2 = 1\}\right| = \max\{\aleph_0, |G|, |H|\}.$$

*Proof.* The cardinality of the free product itself is $\max\{\aleph_0, |G|, |H|\}$, which is therefore an upper bound. The elements of the form $(gh)^n$ are all distinct in $G * H$ for all choices $g \in G \smallsetminus \{1\}$, $h \in H \smallsetminus \{1\}$, and $n \geq 1$. Since none of them can be involutions, this proves the first claim. For the second part, assume $a \in H \smallsetminus \{1\}$ with $a^2 = 1$. Then the elements $a^{(gh)^n}$ derived by conjugating $a$ with $(gh)^n$ are all distinct and $(a^{(gh)^n})^2 = 1$. □

*Proof of Theorem 10.* Suppose $\lambda < \aleph_0$, then by (2) in Lemma 17 we may assume no $G_n$ contains an involution. Then neither does any group $G_{2n-1} * G_{2n}$. Now $G_{2n-1} * G_{2n}$ and $*_{\kappa_{2n-1}} \mathbb{Z} * *_{\kappa_{2n}} \mathbb{Z}$ have the same cardinality $\max\{\aleph_0, \kappa_{2n-1}, \kappa_{2n}\}$. Hence it is possible, to define inverse preserving bijections

$$\varphi_n : G_{2n-1} * G_{2n} \to *_{\kappa_{2n-1}} \mathbb{Z} * *_{\kappa_{2n}} \mathbb{Z},$$

simply by setting $\varphi(1) := 1$ and taking each pair $x, x^{-1}$ of a nontrivial element together with its inverse to some pair $y, y^{-1}$. And so, by applying Corollary 16 to these maps and using property (3) twice, we see

$$\mathscr{A}(G_n) \simeq \mathscr{A}(G_{2n-1} * G_{2n}) \simeq \mathscr{A}(*_{\kappa_{2n-1}} \mathbb{Z} * *_{\kappa_{2n}} \mathbb{Z}) \simeq \mathscr{A}(*_{\kappa_n} \mathbb{Z}).$$

Now suppose $\lambda = \aleph_0$, then by (1) in Lemma 17 the ordering of the groups can be so arranged that every $G_{2n}$ contains an involution. Then by Lemma 20 the cardinalities of both, the set of involutions and non-involutions in $G_{2n-1} * G_{2n}$, are $\max\{\aleph_0, \kappa_{2n-1}, \kappa_{2n}\}$. The same holds for the group $*_{\kappa_{2n-1}} \mathbb{Z}_2 * *_{\kappa_{2n}} \mathbb{Z}_2$. Between these we can again define an identity-preserving map $\varphi_n$, this time by mapping involutions to involutions and pairs of non-involutions $x \neq x^{-1}$ to pairs of non-involutions. By the same argument as in the other case, we get

$$\mathscr{A}(G_n) \simeq \mathscr{A}(G_{2n-1} * G_{2n}) \simeq \mathscr{A}(*_{\kappa_{2n-1}} \mathbb{Z}_2 * *_{\kappa_{2n}} \mathbb{Z}_2) \simeq \mathscr{A}(*_{\kappa_n} \mathbb{Z}),$$

as stated. □

*Proof of Theorem A.* This can be reduced to the case where each group is countably infinite by taking the free product of consecutive pairs of groups $G_n$, as by (3) in Lemma 17 we may write $\mathscr{A}(G_n) = \mathscr{A}(G_{2n-1} * G_{2n})$.

Then applying Theorem 10 twice shows

$$\circledast(G_{2n-1} * G_{2n}) \simeq \mathscr{A}(*_{\aleph_0} \mathbb{Z}) \simeq \mathscr{A}(\mathbb{Z}),$$

if only finitely many groups $G_n$ contain involutions. Similarly with $\mathbb{Z}_2$ replacing $\mathbb{Z}$, if infinitely many groups have 2-torsion. □

*Proof of Theorem 9.* First note that if $\kappa = 1$, then $\mathscr{A}(G_n)$ is the trivial group. Otherwise we calculate the cardinality of $\mathscr{A}(G_n)$ as follows. The topologist's product consists of words whose domain can be any countable order type. For a fixed order type, there will be at most $(\aleph_0 \cdot \kappa)^{\aleph_0}$ words, mirroring the number of functions from a countable set $L$ to $\bigcup_n G_n$, if each of the countably many groups has order $\kappa$. The cardinality of the set of countable order types is $\mathfrak{c} = 2^{\aleph_0}$, thus $|\circledast_n G_n| \leq \mathfrak{c} \cdot \kappa^{\aleph_0} = \kappa^{\aleph_0}$, and that is also an upper bound for the cardinality of the archipelago group.

For simplicity of notation, assume all $G_n$ to be equal to $G$, and let $g^{(n)}$ denote the instance of a $g \in G$ within $G_n$. We will define an injective set function (not a homomorphism) from $G^{\mathbb{N}}$ into $\mathscr{A}(G_n)$, thus providing a lower bound, $|G^{\mathbb{N}}| = \kappa^{\aleph_0} \leq |\mathscr{A}(G_n)|$. For a sequence $(g_j)_j \in G^{\mathbb{N}}$ of elements in $G$ define

$$\varepsilon((g_j)_j) := g_1^{(1)} g_1^{(2)} g_2^{(3)} g_1^{(4)} g_2^{(5)} g_3^{(6)} g_1^{(7)} g_2^{(8)} g_3^{(9)} g_4^{(10)} \cdots,$$

so to each element of $G^{\mathbb{N}}$ we associate an infinite word of order type $\omega$ in $\circledast_n G_n$. Now if $\varepsilon((g_j)_j)$ and $\varepsilon((h_j)_j)^{-1}$ are congruent modulo $*_n G_n$ that implies $g_j = h_j$ for all $j \in \mathbb{N}$. Therefore $\varepsilon$ composed with the quotient map from $\circledast_n G_n$ to $\mathscr{A}(G_n)$ is injective, and hence $|\mathscr{A}(G_n)| = \kappa^{\aleph_0}$.

Recall that $\beth_0 := \aleph_0$ is the cardinality of the integers and $\beth_{\alpha+1} := 2^{\beth_\alpha}$ is the cardinality of the power set of $\beth_\alpha$. Thus for successor cardinals,

$$(\beth_{\alpha+1})^{\aleph_0} = (2^{\beth_\alpha})^{\aleph_0} = 2^{\beth_\alpha} = \beth_{\alpha+1}.$$

A free group in $\kappa \geq \aleph_0$ many generators has cardinality $\kappa$, thus the free groups generated by $\beth_{\alpha+1}$ many generators are all distinct for different $\alpha \geq 0$. $\square$

Next we prove the mapping properties of $\mathscr{A}(\mathbb{Z})$ stated in Theorem B.

**Proposition 21.** *The group of rational numbers $\mathbb{Q}$ is contained as a subgroup in $\mathscr{A}(\mathbb{Z})$.*

*Proof.* We will content ourselves with the basic idea of finding a subgroup isomorphic to $\mathbb{Q}$, a more general construction can be found in [Hoj13]. Consider as an element in $\circledast_n \mathbb{Z}$ the infinite word

$$w := a_1(a_2(a_3(a_4(\ldots)^5)^4)^3)^2,$$

then modulo the free product $*_n \mathbb{Z}$, one can remove the symbol $a_1$ from the word representation and $w \sim (a_2(a_3(a_4(\ldots)^5)^4)^3)^2 =: w_2^2$, so it is a square. Similarly, $w_2 \sim (a_3(a_4(\ldots)^5)^4)^3 =: w_3^3$, and thus $w \sim w_3^6$ is a sixth power. Proceeding in this manner, $w \sim w_n^{n!}$, so $w$ is a divisible element in $\mathscr{A}(\mathbb{Z})$. Thus it is possible to define a homomorphism $\varepsilon : \mathbb{Q} \to \mathscr{A}(\mathbb{Z})$ by setting $\varepsilon(1) := w$, $\varepsilon(1/2) := w_2$, etc. and extending the map to multiples and (additive) inverses of these. $\square$

*Remark.* Proposition 21 also immediately implies that $\mathscr{A}(\mathbb{Z})$ is not $\aleph_1$-free (i.e. not every countable subgroup is free), so (3) in Theorem B cannot be strengthened in that respect.

*Proof of Proposition 7.* For each $G_n$ choose a connected CW-complex $X_n$ whose fundamental group is $G_n$. We can easily assume that every space $X_n$ has a countable local basis for its topology at its basepoint. By Theorem 5, $\mathscr{A}(G_n)$ is the fundamental group of the space $C_f$ where $f$ is the canonical map from the wedge to the shrinking wedge of the $X_n$, as described in the paragraph preceding Theorem 5. The homotopy type of $C_f$ remains unchanged if we attach a cylinder $(\bigodot_n X_n) \times I$ along the shrinking wedge. Let $x$ denote the wedge point of the shrinking wedge in the new base of the thus modified space $D$ and $G := \pi_1(D, x) \simeq \mathscr{A}(G_n)$ its fundamental group. By construction, $D$ now has a countable neighbourhood basis at $x$. Let $\varphi \in \mathrm{Hom}(G, \mathbb{Z})$. By [CC06, Theorem 4.4(3)], there is a neighbourhood $U$ of $x$ so that the homotopy class of any loop in $U$ is in the kernel of $\varphi$. The space $D$ has the property that every loop in $D$ is homotopic into every neighbourhood of $x$. Thus $\varphi$ is trivial. $\square$

Before demonstrating the torsion-freeness of the archipelago groups, it is helpful to characterize the elements of finite order in the topologist's product.

**Lemma 22.** *If $g \in \circledast_n G_n$ has finite order $k$, then $g$ is conjugate to an element $f$ in some $G_i$.*

*Proof.* This relies on the fact that a torsion element in a free product $A * B$ is conjugate to a torsion element in either $A$ or $B$. Thus using the free decompositions

$$\circledast_n G_n = G_1 * \ldots * G_j * \circledast_{n>j} G_n$$

for $j \in \mathbb{N}$, $g$ is either conjugate to an element in some $G_i$, or can be represented as $g = c_j^{-1} f_j c_j$ with $f_j \in \circledast_{n>j} G_n$ and $c_j \in \circledast_n G_n$, for all $j$. But in the latter case, the projection of $g$ into $*_{n=1}^{j} G_n$ is trivial for all $j$. Since the topologist's product embeds in the inverse limit of these groups, that compels $g$ to be trivial. □

*Proof of Proposition 8.* Consider a torsion element $h \in \mathscr{A}(G_n)$. As this group is a direct limit, we know that for some index $j$ there exists a torsion element $g \in \circledast_{n \geq j} G_n$ representing $h = [g]$ in $\mathscr{A}(G_n)$ with respect to the equivalence relation induced by the direct limit. By the preceding lemma, $g$ is conjugate to an element $f$ in $G_i$ (for some $i \geq j$) by $c \in \circledast_{n \geq j} G_n$. Thus

$$h = [g] = [c^{-1} f c] = [c^{-1}][f][c] = [c^{-1}][c] = 1,$$

so the archipelago group is torsion-free. □

In the case where there is no element of order 2 in the individual groups $G_n$, we can strengthen that, so the archipelago group is not only torsion-free, but also locally free.

**Proposition 23.** *$\mathscr{A}(\mathbb{Z})$ is locally free.*

*Proof.* In this case the groups $G_n = \mathbb{Z}$, so we have the direct limit

$$\mathscr{A}(\mathbb{Z}) = \varinjlim_k \circledast_{i \geq k} \mathbb{Z} = \varinjlim_k \mathrm{HE}_k$$

where $\mathrm{HE}_k = \circledast_{i \geq k} \mathbb{Z}$ is the Hawaiian earring group, using only loops labelled $i \geq k$. Recall that the group HE, being a subgroup of the inverse limit of free groups, is locally free, following the result in [CF59].

So $\mathscr{A}(\mathbb{Z})$ is a direct limit of locally free groups, and thus locally free itself. □

For completeness, we append a proof.

**Lemma 24.** *The direct limit of locally free groups is locally free.*

*Proof.* Let $(G_i, i \in I, \varphi_{ij})$ be a directed system of locally free groups. As a set, the direct limit group can be represented as a quotient of the disjoint union of the $G_i$, $G := \varinjlim_i G_i = \coprod_i G_i / \sim$, where $g_i \in G_i$ and $g_j \in G_j$ are equivalent if $\exists k \geq i, j$ such that $\varphi_{ik}(g_i) = \varphi_{jk}(g_j)$.

Let $X$ be a finite subset of $G$. We may assume $X$ to be a minimal generating set of $\langle X \rangle \leq G$. Then for all $x \in X$ there is $j_x \in I$ and $y_x \in G_{j_x}$, such that $\varphi_{j_x}(y_x) = x$. Take $n$ to be an upper bound of the $j_x$. Then $K := \langle \varphi_{j_x, n}(y_x) : x \in X \rangle$ is a free subgroup of the locally free group $G_n$, and $\varphi_n(K) = \langle X \rangle$. As the set $X$ was chosen to be minimal, for each $k \geq n$, the group $\varphi_{nk}(K)$ has the same rank as $K$, and $\varphi_{nk}$ restricted to $K$ is in fact an isomorphism between free groups. Therefore, also $\langle X \rangle$ is free. □

(Gregory R. Conner) DEPARTMENT OF MATHEMATICS, BRIGHAM YOUNG UNIVERSITY, PROVO, UT, 84602, USA

*E-mail address*: `conner@math.byu.edu`

(Wolfram Hojka) INSTITUTE FOR ANALYSIS AND SCIENTIFIC COMPUTATION, TECHNISCHE UNIVERSITÄT WIEN, VIENNA, AUSTRIA

*E-mail address*: `w.hojka@gmail.com`

(Mark Meilstrup) MATHEMATICS DEPARTMENT, SOUTHERN UTAH UNIVERSITY, CEDAR CITY, UT, 84720, USA

*E-mail address*: `mark.meilstrup@gmail.com`